\documentclass[12pt,reqno]{amsart}
\usepackage{
    amsmath,  amsfonts, amssymb,  amsthm,   amscd, 
    gensymb,  graphicx, comment,  etoolbox, url,
    booktabs, stackrel, mathtools,enumitem, mathdots,  microtype, lmodern,    mathrsfs, graphicx, tikz,  longtable,tabularx, float, tikz, pst-node, tikz-cd, multirow, tabularx, amscd,  bm, array, makecell, diagbox, booktabs,ragged2e, caption, subcaption }
\usepackage{makecell,slashbox}
\usepackage{xcolor}
\usepackage[all]{xy}
\usepackage[pagebackref=true, colorlinks=true, linkcolor=blue, citecolor=blue, urlcolor=blue, breaklinks=true]{hyperref}
\usepackage[utf8]{inputenc}
\usepackage{microtype, fullpage, wrapfig,textcomp,mathrsfs,csquotes,fbb}
\usepackage[colorlinks=true, linkcolor=blue, citecolor=blue, urlcolor=blue, breaklinks=true]{hyperref}
 
\usepackage[capitalise]{cleveref}
\usepackage{todonotes}
\usetikzlibrary{positioning}
\usetikzlibrary{shapes,arrows.meta,calc}
\usetikzlibrary{arrows}

\newtheorem{theorem}{Theorem}[section]

\newtheorem{lemma}[theorem]{Lemma}

\newtheorem{proposition}[theorem]{Proposition}
 
\newtheorem*{theorem*}{Theorem}
\newtheorem*{conjecture*}{Conjecture}
\theoremstyle{definition}
\newtheorem{definition}[theorem]{Definition}

\newtheorem{remark}[theorem]{Remark}

\newtheorem*{notation*}{Notation}
 
\newcommand{\m}{\mathfrak{m}}

\newcommand{\n}{\mathfrak{n}}

\newcommand{\grade}{\operatorname{grade}}
\newcommand{\injdim}{\operatorname{injdim}}

\newcommand{\ann}{\operatorname{ann}}

\newcommand{\Ass}{\operatorname{Ass}}
\newcommand{\height}{\operatorname{height}}

\newcommand{\Supp}{\operatorname{Supp}}

\newcommand{\cd}{\operatorname{cd}}

\newcommand{\Spec}{\operatorname{Spec}}

\newcommand{\coker}{\operatorname{coker}}
\newcommand{\Char}{\operatorname{char}}

\newcommand{\Ext}{\operatorname{Ext}}

\makeatletter
\def\@tocline#1#2#3#4#5#6#7{\relax
  \ifnum #1>\c@tocdepth 
  \else
    \par \addpenalty\@secpenalty\addvspace{#2}%
    \begingroup \hyphenpenalty\@M
    \@ifempty{#4}{%
      \@tempdima\csname r@tocindent\number#1\endcsname\relax
    }{%
      \@tempdima#4\relax
    }%
    \parindent\z@ \leftskip#3\relax \advance\leftskip\@tempdima\relax
    \rightskip\@pnumwidth plus4em \parfillskip-\@pnumwidth
    #5\leavevmode\hskip-\@tempdima
      \ifcase #1
       \or\or \hskip 1em \or \hskip 2em \else \hskip 3em \fi%
      #6\nobreak\relax
    \hfill\hbox to\@pnumwidth{\@tocpagenum{#7}}\par
    \nobreak
    \endgroup
  \fi}
\makeatother

\newcolumntype{x}[1]{>{\centering\arraybackslash}p{#1}}

\begin{document}
\title[]{On a generalization of result due to Peskine and Szpiro-II}
\author[Sayed Sadiqul Islam]{Sayed Sadiqul Islam}
\address{Indian Institute of Technology Bombay, Mumbai, India}
\email{22d0786@iitb.ac.in}
\author{Tony J. Puthenpurakal}
\address{Indian Institute of Technology Bombay, Mumbai, India}
\email{tputhen@gmail.com}

\thanks{}
 
\begin{abstract} 
 Let $(R,\m)$ be a regular local ring and let $I$ be a Cohen-Macaulay ideal of $R$. When $\Char R=p>0$, a result of Peskine and Szpiro proves that $H^i_I(R)=0$ for $i\neq\height I$. This result is not true when $\Char R=0$. In this article, we prove that the Bass numbers of $H^{\height I}_I(R)$ completely determines the vanishing of $H^i_I(R)$ for $i\neq \height I$. This result was previously established by the second author in the equicharacteristic case and more generally for Gorenstein local rings by Hellus and Schenzel. In contrast to the spectral sequence methods employed by Hellus and Schenzel, our proof is elementary and self-contained.
\end{abstract}
\maketitle
\section{Introduction}
Let $(R,\m)$ be a regular local ring. We say that an ideal $I$ of $R$ is Cohen-Macaulay if $R/I$ is Cohen-Macaulay local ring. We recall the definition of Peskine-Szpiro ideal.
\begin{definition}
    An ideal $I$ of $R$ is said to be a Peskine-Szpiro ideal of $R$ if the following hold.
    \begin{enumerate}
        \item $I$ is a Cohen-Macaulay ideal of $R$.
        \item $H^i_I(R)=0$ for all $i\neq \height I$.
    \end{enumerate}
\end{definition}
Note that as $\height I=\grade I$, we have $H^i_I(R)=0$ for $i<\height I$. Thus a Cohen-Macaulay ideal $I$ is Peskine-Szpiro ideal if $H^i_I(R)=0$ for $i>\height I$. The next result is due to Peskine and Szpiro \cite[Proposition III.4.1]{PS73}.
\begin{theorem}\label{Peskineszpiro result}
    Let $R$ be a regular domain of characteristic $p>0$. If $I$ is a Cohen-Macaulay ideal of $R$, then $H^i_I(R)=0$ for all $i\neq \height I$.
\end{theorem}
The above result does not hold over regular rings of characteristic zero, for instance, see \cite[Example 21.31]{ILL+07}. This motivated the second author to prove the following surprising result.
\begin{theorem} \cite[Theorem 1.2]{Put24}\label{put24 main result}
    Let $(R,\m)$ be a regular local ring of dimension $d$ containing a field $K$. Let $I$ be a Cohen-Macaulay ideal of height $g$. The following conditions are equivalent.
    \begin{enumerate}[label=(\alph*), font=\upshape]
        \item I is a Peskine-Szpiro ideal of $R$.
        \item If $P\in V(I)$,
         $$
\mu_i(P,H^g_I(R)) =
\begin{cases}
1 & \text{if } i=\height P-g,\\
0, & \text{otherwise }.
\end{cases}
$$
    \end{enumerate}
\end{theorem}
The proof of Theorem \ref{put24 main result} uses the fact that $\injdim_R H^i_I(R)\leq \dim\Supp H^i_I(R)$ which was prove by Huneke and Sharp \cite{HS93} when $\Char K=p>0$ and by Lyubeznik \cite{Lyu93} when $\Char K=0$. Unfortunately, when the regular ring  does not contain a field, the bound on injective dimension of local cohomology modules are different. In \cite{CZ98}, Zhou proved that  $\injdim_RH^i_I(R)\leq \dim_R(\Supp_R H^i_I(R))+1$ when $R$ is an unramified regular local of mixed characteristic and $I$ is an ideal of $R$. The aim of this article is to establish Theorem \ref{put24 main result} without the assumption that $R$ contains field. More precisely, we prove the following.
\begin{theorem}( Theorem \ref{proof of main theorem} below)
    Let $(R,\m)$ be a regular local ring of dimension $d$. Let $I$ be a Cohen-Macaulay ideal of height $g$. The following conditions are equivalent.
    \begin{enumerate}[label=(\alph*), font=\upshape]
        \item I is a Peskine-Szpiro ideal of $R$.
        \item If $P\in V(I)$,
         $$
\mu_i(P,H^g_I(R)) =
\begin{cases}
1 & \text{if } i=\height P-g,\\
0, & \text{otherwise }.
\end{cases}
$$
    \end{enumerate}
\end{theorem}
We remark that the above 
result follows as a special case of a more general result due to Hellus and Schenzel \cite[Theorem 3.2]{HS08}. However, our proof is elementary and avoids a spectral sequence argument.
\section{Cohen-Macaulay ideals of small dimensions}
Let $R$ be a Noetherian ring and let $I$ be an ideal of $R$. The local cohomological dimension of $I$, denoted as $\cd(I)$ is defined as follows:
$$\cd(I)=\sup\{j\mid H^j_I(R)\neq 0\}.$$
Let us first recall the first vanishing theorem of local cohomology which is commonly known as Hartshorne-Lichtenbaun Vanishing Theorem (HLVT for short) \cite{Har68}.
\begin{theorem}(Hartshorne-Lichtenbaum vanishing theorem)\label{HLVT}
    Let $(R,\m)$ be a complete local domain. If $I$ is an ideal of $R$ which is not $\m$-primary, then $\operatorname{cd}(I)\leq \dim R-1$.
\end{theorem}
Let us now recall the second vanishing theorem of local cohomology modules.
\begin{theorem}
    Suppose $(R,\m)$ is a complete regular local ring of dimension $d$ with separably closed residue field. Then, $$\cd(R/I)\leq d-2 \iff \dim(R/I)\geq 2 \ \text{and} \Spec(R/I)\setminus\{\m\}\ \text{is connected}.$$
\end{theorem}
The equal characteristic zero and positive characteristic cases were established by Ogus \cite{Ogu73} and Peskine–Szpiro \cite{PS73}, respectively. Later, Huneke and Lyubeznik \cite{HL90} gave a characteristic-free proof of the second vanishing result. The unramified mixed characteristic case was proved by Zhang \cite{Zha25} and the remaining ramified mixed characteristic case was recently settled by Scheffelin \cite{Sch26}. We use the second vanishing theorem to prove the following proposition.
\begin{proposition}\label{PS ideals of small dimension}
    Let $(R,\m)$ be a regular local ring of dimension $d$ and let $I$ be a Cohen-Macaulay ideal of $R$ of height $g$ with $d-g\leq 2$. Then $I$ is a Peskine-Szpiro ideal.
\end{proposition}
\begin{proof}
    When $\Char R=p>0$, the result follows from \cite[Proposition III.4.1]{PS73}. Therefore we assume that $\Char R=0$.

    If $g=d$, then $I$ is $\m$-primary and therefore by Grothendieck vanishing theorem $H^i_I(R)=0$ for $i> \height I$. So, $I$ is a Peskine-Szpiro ideal.

    Now consider $g=d-1$. Therefore $I$ is not  $\m$-primary and hence $I\hat{R}$ is not $\m\hat{R}$-primary. By Hartshorne-Lichtenbaum theorem, $H^d_{I\hat{R}}(\hat{R})=0$. By faithful flatness we get $H^d_I(R)=0$.

    Now let $g=d-2$. Choose a faithfully flat extension $(S,\n)$ of $R$ such that $\m S=\n$, $S$ is complete and $S/\n$ is algebraically closed. We note that $(S,\n)$ is regular local, $S/I S$ is Cohen-macaulay and $\dim S/I S=2$. Therefore $\Spec(S/IS)\setminus\{\n\}$ is connected by \cite[Theorem 15.7]{ILL+07}. By the second vanishing theorem we get  $H^{d-1}_{IS}(S)=0$ and this gives $H^{d-1}_I(R)=0$. By a similar argument to the previous case we also get $H^d_I(R)=0$.
    \end{proof}
    The proof of the following result is identical to that of \cite[Proposition 2.1]{Put24} and is omitted to avoid repetition.
\begin{proposition}\label{localization and quotient is PS}
    Let $(R,\m)$ be a regular local ring. Let $I$ be a Peskine-Szpiro ideal of $R$ of height $g$. Then, the following statements are true.
    \begin{enumerate}[label=(\alph*), font=\upshape]
        \item Let $P$ be an ideal in $R$ containing $I$. Then $I_P$ is a Peskine-Szpiro ideal of $R_P$.
        \item Assume $g<\dim R$. Let $x\in \m-\m^2$ be $R/I$-regular. Then the ideal $I+(x)/(x)$ is a Peskine-Szpiro ideal of $R/(x)$.
    \end{enumerate}
\end{proposition}
\section{Proof of the main theorem}
The following fact is well known, but we include a proof for completeness, as we are unable to find a suitable reference.
\begin{lemma}\label{easy lemma}
    Let $(R,\m)$ be a Noetherian local ring and let $M$ be an $R$-module. Assume that $\Supp M=\{\m\}$ and $x\in \m$. Then $\ker(M\xrightarrow{x}M)\neq 0$.
\end{lemma}
\begin{proof}
    Let $z\in M$ be non-zero. Then $\Supp(Rz)=V(\ann(Rz))=\{\m\}$. This implies $\sqrt{\ann(Rz)}=\m$ and hence $\m^tz=0$ for some $t\geq 1$. In particular, $x^tz= 0$. Choose a smallest $t$ with this property. Then $y=x^{t-1}z\neq 0$ and $y\in \ker(M\xrightarrow{x}M)$. 
\end{proof}
We also note the following remark.
\begin{remark}\label{remark about ass prime}
   Let $R$ be a regular ring and let $J$ be an ideal of height $g$ of $R$. Then $\Ass H^g_J(R)=\{P\mid P\supseteq J\ \text{and}\  \height P=g\}$. Also for any such prime $P$ we have $\mu_0(P,H^g_J(R))=1$ since $H^g_J(R)_P=H^g_{JR_P}(R_P)=H^g_{PR_P}(R_P)=E_{R_P}(R_P/PR_P)$.
\end{remark}
Let us also recall the following result due to Rees \cite[Lemma 3.1.16]{BH93}.
\begin{lemma}\label{rees lemma}
    Let $S$ be a commutative ring and let $M$ and $N$ be $S$-modules (not necessarily finitely generated). Let $x\in S$ be $S\oplus M$ regular with $xN=0$. Then $$\Ext^i_S(M,N)\cong \Ext^i_{S/xS}(M,N/xN).$$
\end{lemma}
We are now in a position to prove the main theorem. While the direction $(a) \implies (b)$ is analogous to \cite[Theorem 1.2]{Put24}, we provide the complete proof here to ensure the paper is self-contained. Note that the proof of the converse direction $(b) \implies (a)$ does not follow from \cite[Theorem 1.2]{Put24} and requires modification.
\begin{theorem} \label{proof of main theorem}
    Let $(R,\m)$ be a regular local ring of dimension $d$. Let $I$ be a Cohen-Macaulay ideal of height $g$. The following conditions are equivalent.
    \begin{enumerate}[label=(\alph*), font=\upshape]
        \item I is a Peskine-Szpiro ideal of $R$.
        \item If $P\in V(I)$,
         $$
\mu_i\left(P,H^g_I(R)\right) =
\begin{cases}
1 & \text{if } i=\height P-g,\\
0, & \text{otherwise }.
\end{cases}
$$
    \end{enumerate}
\end{theorem}
\begin{proof}
    We first prove $(a)\implies (b)$. Let $I$ be a Peskine-Szpiro ideal. We prove the result by induction on $d-g$. If $d-g=0$, then $\sqrt{I}=\m$ and hence $H^d_I(R)=H^d_\m(R)=E_R(R/\m).$

    Now assume $d-g=1$. Let $P\in V(I)$ with $\height P=d-1$. By Remark \ref{remark about ass prime}, $\mu_0(P,H^{d-1}_I(R))=1$ and $\mu_i(P,H^{d-1}_I(R))=0$ for $i\geq 1$. Now let $P=\m$. Then by Remark \ref{remark about ass prime}, $\mu_0(\m,H^{d-1}_I(R))=0$. Choose $x\in \m\setminus \m^2$ which is $R/I$-regular. Let $\overline{R}=R/(x), \n=\m/(x) \ \text{and} \ J=I\overline{R}=I+(x)/(x)$. Note that $J$ is $\n$-primary. The short exact sequence $0\to R\xrightarrow{x}R\to \overline{R}\to 0$ induces the following long exact sequence on the local cohomology
    \begin{equation*}
        0\to H^{d-1}_I(R)\xrightarrow{x}H^{d-1}_I(R)\to H^{d-1}_J(\overline{R})\to 0.
    \end{equation*}
Therefore,
 $$
\mu^R_i\left(\m,H^{d-1}_I(R)\right) = \mu^{\overline{R}}_{i-1}\left(\n,H^{d-1}_J(\overline{R})\right)=
\begin{cases}
1 & \text{if } i=1,\\
0, & \text{otherwise }.
\end{cases}
$$
Thus the result follows in this case.

Now consider the case $d-g\geq 2$. Now let $P\in V(I)$ be such $\height P=g+r$ with $P\neq \m$. Then by Proposition \ref{localization and quotient is PS}, $I_P$ is a Peskine-Szpiro ideal of $R_P$. Note that $\dim R_p-g<d-g$. Therefore by the induction hypothesis 
 $$
\mu_i\left(P,H^{g}_I(R)\right) = \mu_{i}\left(PR_P,H^{g}_{I_P}(R_P)\right)=
\begin{cases}
1 & \text{if } i=r,\\
0, & \text{otherwise }.
\end{cases}
$$
We now consider the case $P=\m$. By Remark \ref{remark about ass prime}, $\mu_0(\m,H^g_I(R))=0$.  Choose $x\in \m\setminus \m^2$ which is $R/I$-regular. Let $\overline{R}=R/(x), \n=\m/(x) \ \text{and} \ J=I\overline{R}=I+(x)/(x)$. By Proposition \ref{localization and quotient is PS}, $J$ is a height $g$ Peskine-Szpiro ideal in $\overline{R}$. The short exact sequence $0\to R\xrightarrow{x}R\to \overline{R}\to 0$ induces the following long exact sequence on the local cohomology
    \begin{equation*}
        0\to H^{g}_I(R)\xrightarrow{x}H^{g}_I(R)\to H^{g}_J(\overline{R})\to 0.
    \end{equation*}
    The result now follows from the induction hypothesis ($\dim \overline{R}-g=d-1-g$) and using the fact $\mu_i(\m,H^g_I(R))= \mu_{i-1}\left(\n,H^{d-1}_J(\overline{R})\right)$.
\medskip

    We now prove $(b)\implies (a)$. By Theorem \ref{Peskineszpiro result}, we may assume $\Char R=0$. We apply induction on $d-g$. If $d-g\leq 2$, the result follows from Proposition \ref{PS ideals of small dimension}. Therefore, we assume that $d-g\geq 3$. Let $P\in V(I)$ with $P\neq \m$. Then, $I_P$ is a Cohen-Macaulay ideal of $R_P$ of height $g$ and satisfies (b). As $\dim R_P-g<d-g$, induction hypothesis implies that $I_P$ is a Peskine-Szpiro ideal of $R_P$. Therefore, $H^i_{I_P}(R_P)=0$ for all $i>g$. It follows that $\Supp_R H^i_I(R)\subseteq\{\m\}$ for all $i>g$.

    Choose $x\in \m\setminus \m^2$ which is $R/I$-regular. Let $\overline{R}=R/(x), \n=\m/(x) \ \text{and} \ J=I\overline{R}=I+(x)/(x)$. The short exact sequence $0\to R\xrightarrow{x}R\to \overline{R}\to 0$ induces the following long exact sequence on the local cohomology
    \begin{equation}\label{les of lc}
        0\to H^g_I(R)\xrightarrow{x}H^g_I(R)\to H^g_J(\overline{R})\to H^{g+1}_I(R)\xrightarrow{x}H^{g+1}_I(R)\to \cdots.
    \end{equation}

Let $V^r:=\coker(H^r_I(R)\xrightarrow{x}H^r_I(R))$ and $L^r:=\ker(H^{r}_I(R)\xrightarrow{x}H^{r}_I(R))$ for all $r$. Therefore we have the following two short exact sequences \begin{equation}\label{equation 1}
    0\to H^g_I(R)\xrightarrow{x}H^g_I(R)\to V^g\to 0
\end{equation}
and \begin{equation}\label{equation 2}
    0\to V^g \to H^g_J(\overline{R})\to L^{g+1} \to 0.
\end{equation}
Note that $\Supp_R L^{g+1}\subseteq \Supp_R H^{g+1}_I(R)\subseteq\{\m\}.$ We first claim that $H^{g+1}_I(R)= 0$.

If possible let $H^{g+1}_I(R)\neq 0$. Then By Lemma \ref{easy lemma}, $\Supp_R H^{g+1}_I(R)=\Supp_R L^{g+1}=\{\m\}$. 
Since $R/I$ is Cohen-Macaulay and  $\overline{R}/J=R/I+(x)=(R/I)/x (R/I)$, we see that $J$ is Cohen-Macaulay ideal of $\overline{R}$. Also $\height J=\dim \overline{R}-\dim \overline{R}/J=\height(I)=g$. By Remark \ref{remark about ass prime}, $\n\notin\Ass H^g_J(\overline{R})$ as $d-1\neq g$ and hence $\mu^{\overline{R}}_0\left(\n,H^g_J(\overline{R})\right)=0$. Therefore, the short exact sequence (\ref{equation 2}) gives $$\mu^{\overline{R}}_1\left(\n,V^g\right)\geq \mu^{\overline{R}}_0(\n,L^{g+1})>0.$$
But $\mu^{\overline{R}}_1(\n,V^g)=\mu^R_2\left(\m,H^g_I(R)\right)$ by Lemma \ref{rees lemma} and therefore by condition (b), $d-g=2$ which is a contradiction. Consequently, $H^{g+1}_I(R)=0$.

Since $H^{g+1}_I(R)=0$, from (\ref{les of lc}) we get the following short exact sequence $$ 0\to H^g_I(R)\xrightarrow{x}H^g_I(R)\to H^g_J(\overline{R})\to 0.$$
Note that $J$ is a Cohen-Macaulay ideal of $\overline{R}$ and satisfies condition (b) of our theorem. Hence by the induction hypothesis, $H^i_J(\overline{R})=0$ for $i>g$. Therefore, from (\ref{les of lc}), we get $L^i=0$ for all $i\geq g+2$. This implies $H^i_I(R)=0$ for $i\geq g+2$ by Lemma \ref{easy lemma} as $\Supp_R H^i_I(R)=\{\m\}$ for $i\geq g+2$.
\end{proof}
\bibliographystyle{alpha}
 
\end{document}